\documentclass[12pt]{article}
\usepackage{amsmath,amssymb,latexsym}
\usepackage{times}
\usepackage{graphicx}

\newtheorem{theorem}{Theorem}[section]
\newtheorem{remark}{Remark}[section]

\numberwithin{equation}{section}

\begin{document}

\title{Word Representations of $m\times n\times p$ Proper Arrays}
\date{\today}
\author{Jocelyn Quaintance\\quaintan@temple.edu}
\maketitle
\begin{abstract}
\noindent Let $m\neq n$.  An $m\times n\times p$ {\it proper array} is a three-dimensional array composed of directed cubes that obeys certain
constraints.  Because of these constraints, the $m\times n\times p$ proper arrays may be classified via a schema in which each $m\times n\times p$ proper array is associated with a particular $m\times
n$ planar face.  By representing each connencted component present in the $m\times n$ planar face with a distinct letter, and the position of each outward pointing connector by a circle, an $m\times n$
array of circled letters is formed.  This $m\times n$ array of circled letters is the {\it word representation} associated with the
$m\times n\times p$ proper array.  The main result of this paper involves the enumeration of all $m\times n$ word representations modulo symmetry, where the symmetry is derived from
the group $D_2 = C_2\times C_2$ acting on the set of word representations.  This enumeration is achieved by forming a linear combination of four exponential generating functions, each of which is derived
from a particular symmetry operation.  This linear combination counts the number of partitions of the set of $m\times n$ words representations that are inequivalent under $D_2$
\end{abstract}
\newpage
\section*{Introduction}
\subsection{History of the Problem}
\noindent This is a continuation of the author's work on the enumeration of three-dimensional proper arrays [6].  The building block of a proper array is a {\bf directed cube}, where a directed cube
is a cube which has five faces containing circular indentations while the sixth face has a cylindrical plug or {\bf connector}
\begin{center}
\includegraphics [width = 3.0cm, height = 3.0cm]{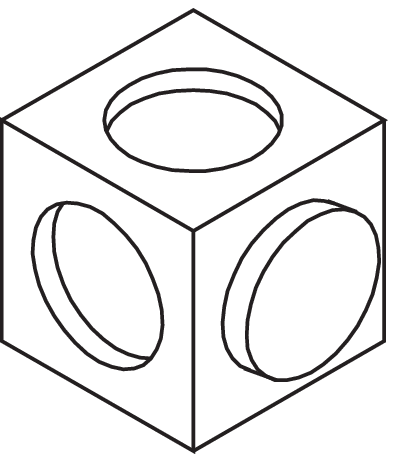}
\end{center}
\begin{center}
\small{ Figure 1: A directed cube.  In this illustration, the connector extends from the right face while the other five faces contain circular indentations.}
\end{center}
\noindent  Three-dimensional arrays can be built by fitting the connector of one directed cube into an indentation of another.  An $m\times n\times p$ array of directed cubes is composed of $p$
layers, with $m$ rows and $n$ columns in each layer.  If the $m\times n\times p$ array of directed cubes obeys certain connectivity constraints, it is said to be an
{\boldmath $m\times n\times p$} {\bf proper array} [6].  An $m\times n\times p$ proper array may be classified by its outermost $m\times n$ layer of directed cubes, which is
referred to as the {\bf preferred face} of the $m\times n\times p$ proper array [6].  The preferred face of the $m\times n\times p$ proper array contains all its path-connected subsets and all its
outward extending connectors.\\\\ 
 \noindent The preferred face does not uniquely determine the number $m\times n\times p$ proper arrays since two different $m\times n\times p$ proper arrays may have the same
preferred face.  However, the preferred face is the state of a transition matrix {\boldmath $M_{m\times n}$} [1],[2],[12].  For a fixed $m$ and $n$ and arbitrary $p$, {\boldmath $M_{m\times n}$} is an
$r\times r$ matrix that describes how an $m\times n\times p$ proper array evolves into an $m\times n\times (p+1)$ proper array.  The number $r$ is the number of possible preferred faces, modulo
symmetry, associated with $m\times n\times p$ proper arrays [1],[2],[12].  Let {\boldmath $S_p$} be the $1\times r$ vector
$(s_{p,i})$, where $s_{p,i}$ counts, modulo symmetry, the number of $m\times n\times p$ proper arrays having the $i$th kind of preferred face.  We find {\boldmath $S_pM_{m\times n}$} $=$ {\boldmath
$S_{p+1}$}.  For example, in the case of $2\times 1\times p$ proper arrays, there are four possible preferred face structures.  Hence, {\boldmath $M_{2\times 1}$} is a
$4\times 4$ matrix.  Let {\boldmath $S_1$} $ = [1,\,\,2, \,\, 0,\,\, 1]$, where the entries of {\boldmath $S_1$} count the number of $2\times 1\times 1$ proper arrays with a particular preferred face 
structure.  Note, there are $1 + 2 + 0 + 1 = 4$ distinct $2\times 1\times 1$ proper arrays.  Then {\boldmath $S_1M_{2\times 1}$} $=$ {\boldmath $S_2$} $ = [2,\,\,8,\,\,0,\,\, 3]$, where the entries of
{\boldmath $S_2$} count the number of $2\times 1\times 2$ proper arrays with a particular preferred face structure.  Note, there are
$2 + 8 + 0 + 3 = 13$ distinct $2\times 1\times 2$ proper arrays.  \\\\ 
\newpage
\noindent For a fixed $m$, fixed $n$, and arbitrary $p$, we construct {\boldmath $M_{m\times n}$} via a computer program.  Thus, one of the main research question associated with the enumeration of
$m\times n\times p$ proper arrays is predicting the size of $r$, where $r$ is the number of preferred faces, modulo symmetry.  We should mention that if $m\neq n$, the symmetry in question consists of
four maps, each of which maps the preferred face onto itself, i.e. the symmetry group is $D_2$.  By having a formula for $r$, we will be able to determine the amount of time and memory necessary for the
computation of {\boldmath $M_{m\times n}$} [5].  The main result of this paper is the first step in determining such a formula.  In particular,  Theorem 2.1 provides a formula, in the form of a linear
combination of four exponential generating functions, for calculating, modulo  $D_2$ symmetry, an upper bound on the size of $r$.  The techniques used to derive the generating
functions are similar in nature to the techniques used by Yoshinaga and Mori [7] and David Branson [8].  The second result provides generating functions which count the number of $m\times n\times p$
proper arrays whose preferred face is fixed by the four symmetry maps of $D_2$.
\subsection{Connections to Percolation Theory}
\noindent Before continuing, we would like to provide a physical interpretation for an $m\times n\times p$ proper array.  If the center of each directed cube is represented as a point of
$Z^3$, with the connectors providing a network of open paths between the various vertices, then the $m\times n\times p$ proper array is a model of bond percolation [10 Section 1], [11
P.16].  Furthermore, this bond percolation has the following property: there exists an open cluster that allows water to flow from the back face of the $m\times n\times p$ array to its front face
[10 P.2].  Moreover, since there is a relationship between bond percolation and the Ising Model [10 P.8], an open question involves the exact nature of the connection between $m\times n\times p$ proper
arrays and Ising Models.  At this point in time, we have not investigated this connection between the Ising Model and the three-dimensional proper arrays and leave it as an open question for future
research. 
\section{Proper Arrays and Word Representations}
\noindent We are now ready to give the definition of an $m\times n\times p$ proper array.\\\\
\noindent {\bf Definition:} A (connected) {\bf component} of an array of directed cubes is a (path) connected subset of the array.\\\\
\noindent {\bf Definition:} Let $A$ represent an $m\times n\times p$ array of directed cubes.  Orient $A$ so that one planar $m\times n$ face has a center at $(0,1,0)$ and is perpendicular to the $y$
axis.  This $m\times n$ planar face is the {\bf preferred face} of $A$.\\\\

\noindent {\bf Definition:} An {\bf $m\times n\times p$ proper array} is a three-dimensional array of directed cubes that obeys the following two conditions:
\begin{itemize}
\item[1.]  The array is {\bf 5-way flat}.  This means only the preferred face has outward pointing connectors.
\begin{center}
\includegraphics[width=1.75in, height=2.0in]{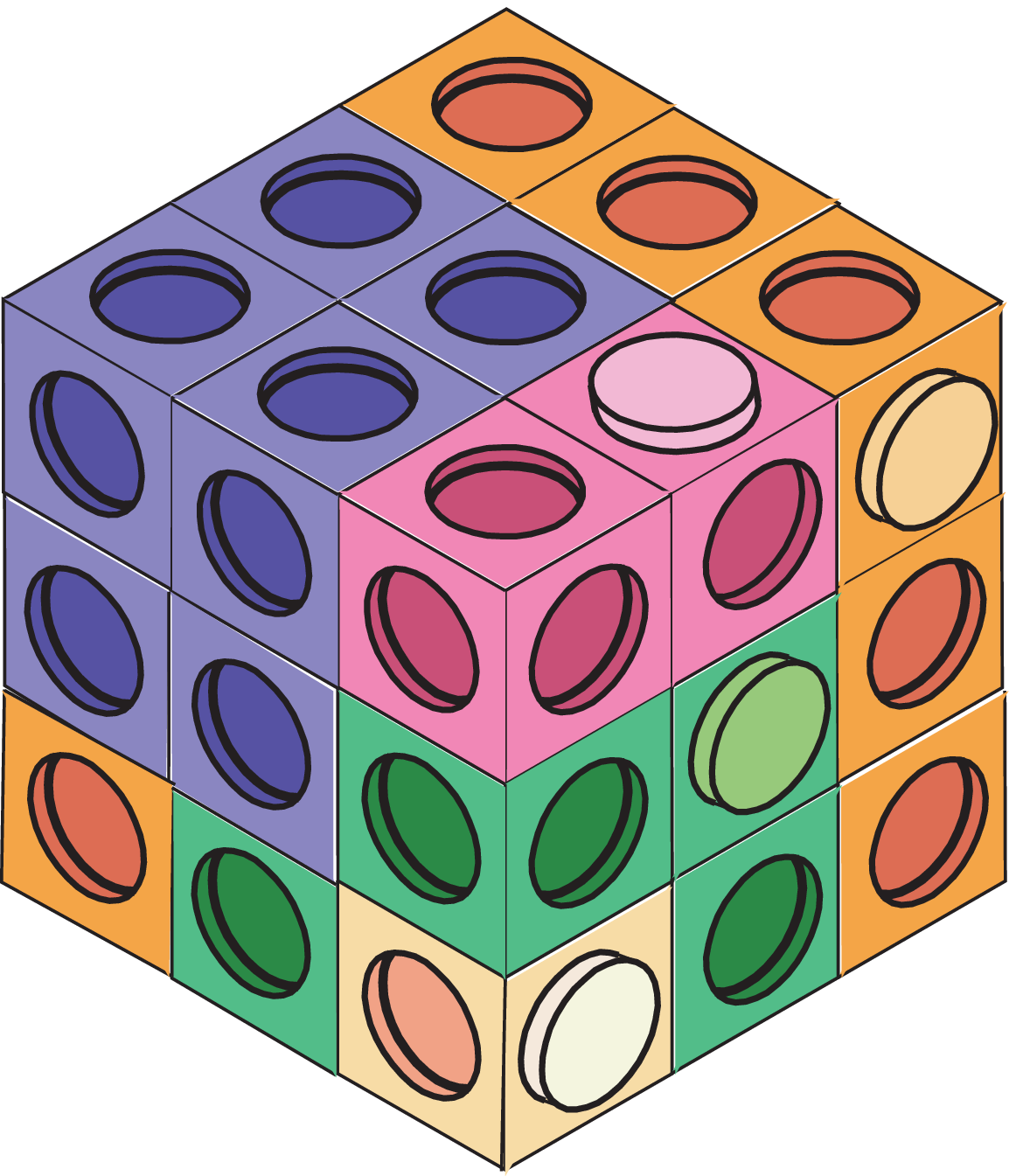}
\end{center}
\begin{center}  
\small{Figure 1.1: This $3\times 3\times 3$ array of directed cubes is not 5-way flat since the pink component has a conncetor that extends from the top face.}
\end{center}
\item[2.]  We define an {\bf island} to be a component that does not reach the preferred face. The array \underline{can not} have islands.
\begin{center}
\includegraphics[width=1.75in, height=2.0in]{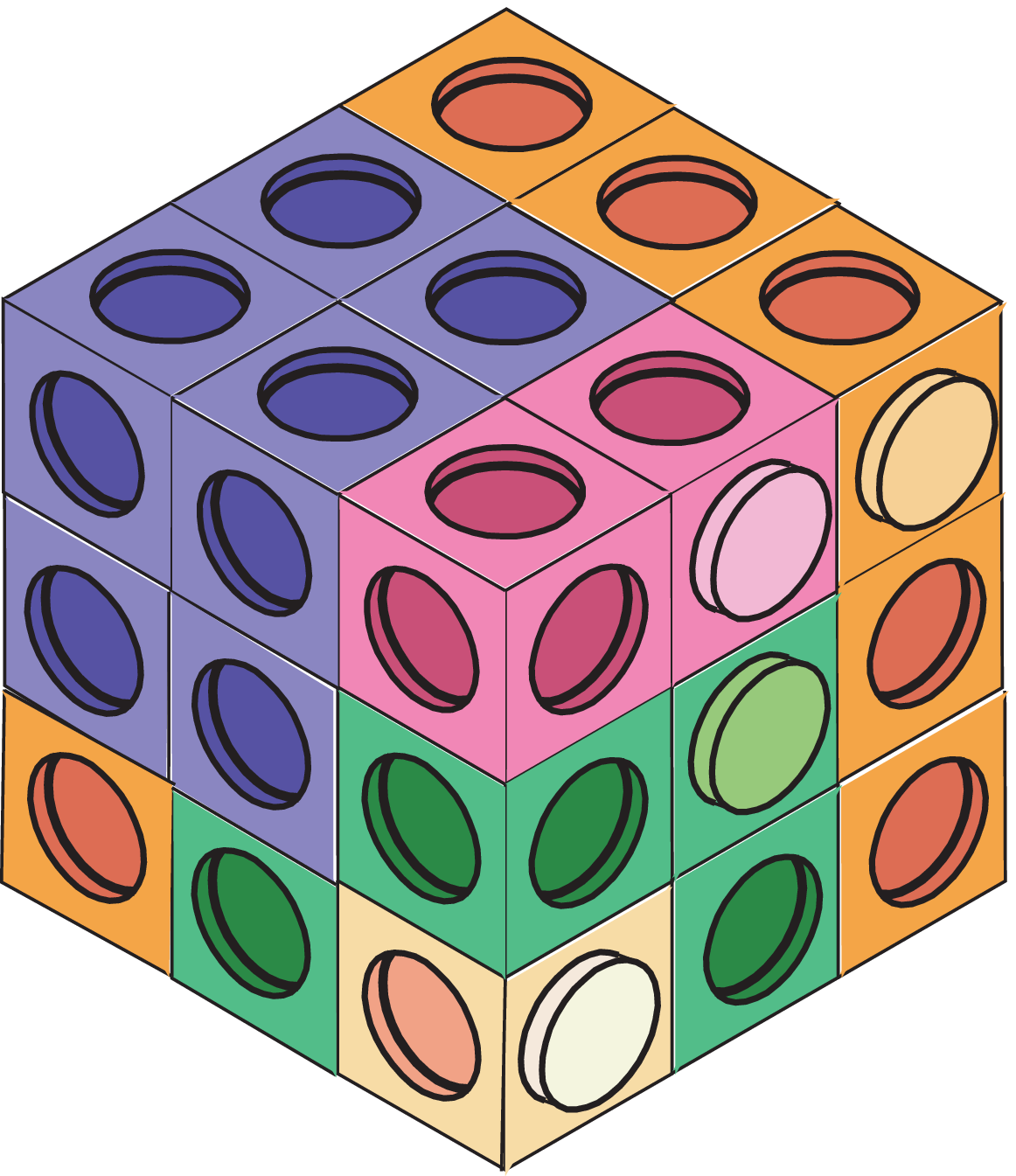}
\end{center}
\begin{center} 
\small{Figure 1.2: The purple component forms an island since it does not reach the preferred face.}
\end{center}
\end{itemize}
\noindent An $m\times n\times p$ proper array may be associated with its preferred face, since the preferred face corresponds to the state of the transition matrix {\boldmath $M_{m\times n}$}.  In
particular, by letting each distinct component be denoted by a distinct letter, and any outward pointing connectors be depicted as circles around the letters, the geometry of the preferred face
is recorded as $m\times n$ array of circled letters.  This $m\times n$ array of circled letters is called the {\boldmath $m\times n$} {\bf word
representation} of the $m\times n\times p$ proper array.  The word representation encodes the numbers of connected components that appear in the $m\times n\times p$ proper array and the outward
extending connectors of $m\times n\times p$ proper array.  
\begin{center}
\includegraphics[width=3.5in, height=2.0in]{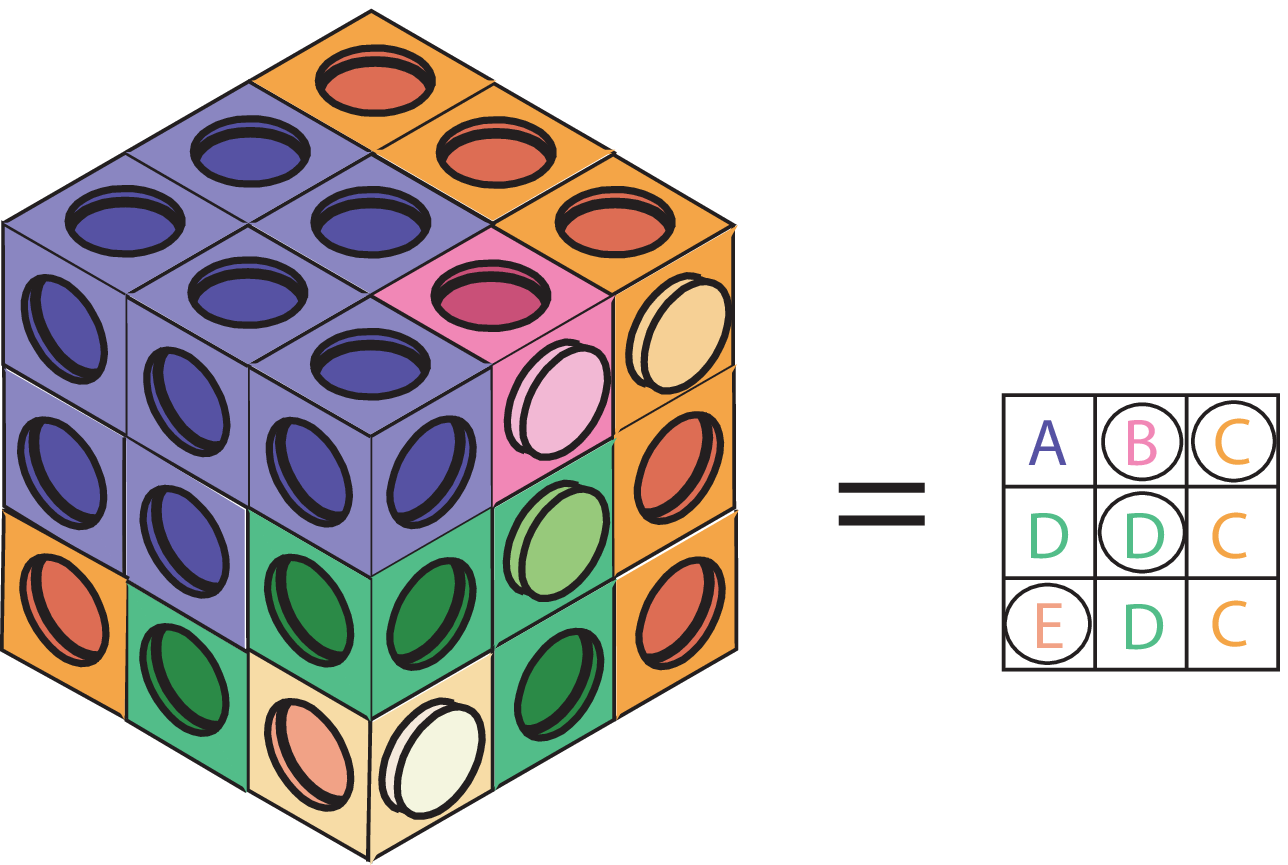}
\end{center}
\begin{center} 
\small{Figure 1.3: An example of a $3\times 3\times 3$ proper array with its associated word representation.  This word representation records
the components as colored letters and the outward pointing connectors as circles.}
\end{center}  
\begin{remark}
Word Representations are partitions of the preferred face.  Hence, two word representations are regarded as the same if they have the same components but their letters are different.  In other words,
the letters are just labels for the components and the labels are unimportant.  For example, take the word representation provided in Figure 1.3.  If we replace the $A$ with an $X$, the resulting
word representation is considered to be the same as the orginal word representation.
\end{remark} 
\noindent By applying path connected arguments to the $m\times n\times p$ proper array, it is fairly straight forward exercise to prove Remark
1.2 ([1],[2]).  Note that Remark 1.2 restricts the positions of the cirlces in the word representation. 
\begin{remark}
Each connected component in the $m\times n\times p$ proper array can have, at most, one outward pointing connector.  In other words, an $m\times n$ word representation obeys the following property: There
is at most one circle present in each collection of letters that represent a single connected component.
\end{remark}
\noindent For the remainder of this paper, we assume all word representations obey Remark 1.2.\\\\
\noindent The definition of an $m\times n\times p$ proper array is motivated by the desire to ensure that a particular $m\times n$ face of the proper array corresponds to the state of the transition
matrix [1],[2],[12].  In order to see how the definition of proper array guarantees this correspondence, we will take a moment to describe how the transition matrix constructs
proper arrays.  Fix $m$ and $n$.  The matrix multiplication {\boldmath $S_{p-1}M_{m\times n}$} corresponds to attaching an $m\times n$ layer of directed cubes onto a previous constructed $m\times n\times
(p-1)$ proper array.  The concept of 5-way flat ensures that this conncection could only occur on the preferred face.  The concept of no islands ensures there are no connected subsets of the
$m\times n\times (p-1)$ proper array that do not reach its preferred face.  Hence, all the connected subsets of the $m\times n\times (p-1)$ proper array will be affected by the attachment of the new
$m\times n$ layer and will percolate throughout the newly constructed $m\times n\times p$ proper array.  Since this percolation occurs, we are able to record the state of construction by simply looking at
the geometry of the $m\times n$ preferred face of the newly build $m\times n\times p$ proper array.  Without the conditions of 5-way flat and no islands, this correspondence between preferred face and
construction state would not exist.
\section{Enumerating Word Representations}
\noindent  For the remainder of this paper, we can ignore the three-dimensional context provided by the proper array and work in the two-dimensional
setting of the word representation.\\\\
\noindent Our goal is to fix $m$ and $n$, and determine, modulo symmetry, an upper bound for the number of word representations
associated with $m\times n\times p$ proper arrays.  If $m\neq n$, the symmetry equivalence is determined by $D_2$.   
 \begin{center}
\includegraphics[width=3.0in, height=3.5in]{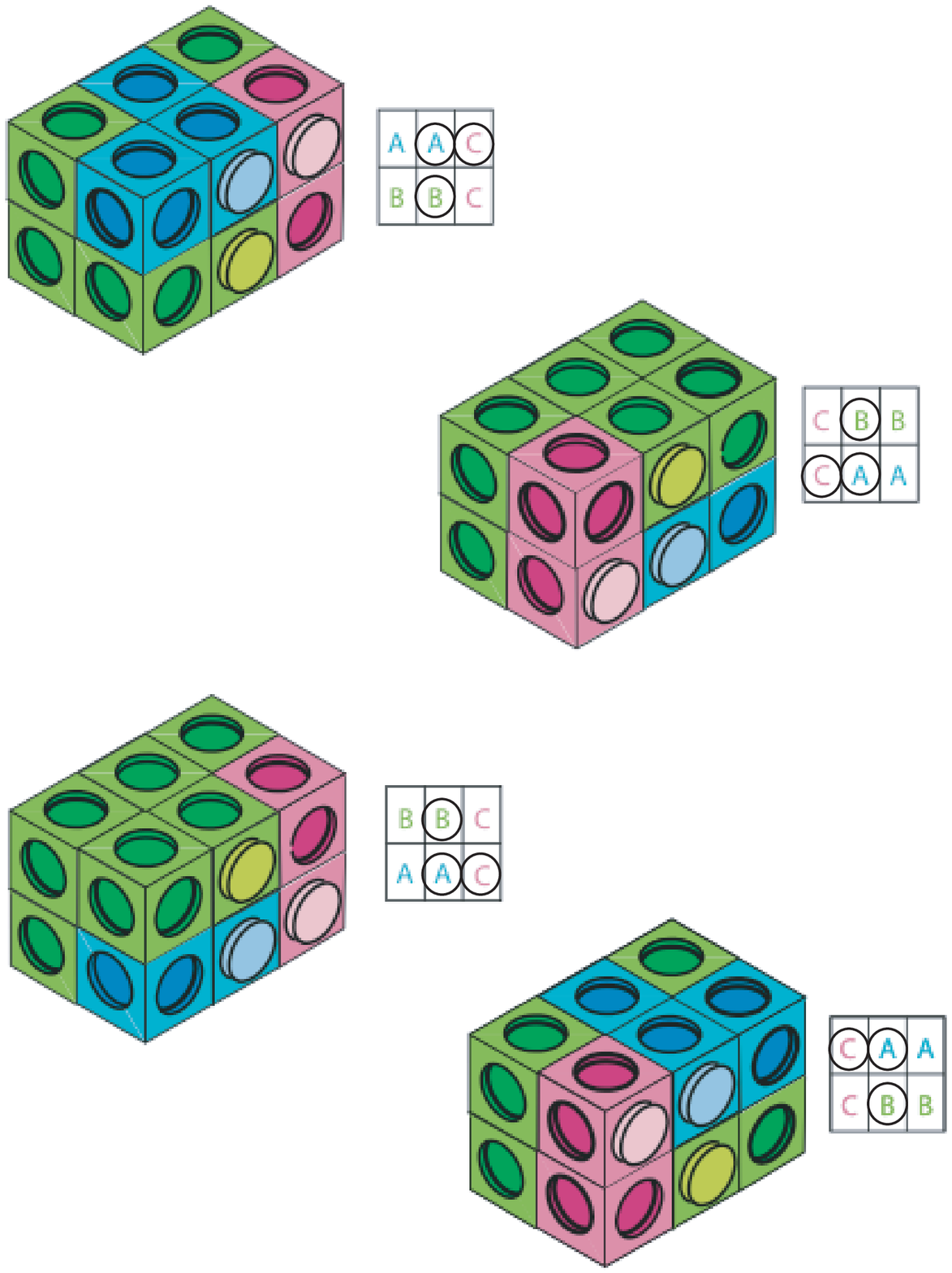}
\end{center}
\begin{center}
\small{Figure 2.1:  The four equivalent versions of a $2\times 3\times 2$ array.  Let the diagram in the upper left corner represent the orginal array.  Then by reading in a clockwise manner, the next
diagram depicits a $180^{\circ}$ rotation of the orginal array.  The bottom right corner depicits a vertical reflection of the orginal array.  The remaining diagram depicits a horizontal reflection of
the orginal array.  These four symmetry images are only counted once in the enumeration procedure.} 
\end{center}
\newpage
\noindent {\bf Definition:}  An {\bf {\boldmath $m\times n$} circled letter array} is an $m\times n$ arrangement of letters and circles which obeys the following two conditions.
\begin{enumerate}  
\setlength{\itemsep}{-\itemsep}
\item[1.] Each distinct letter labels an element in the partition of the $m\times n$ grid of squares.
\item[2.] Each element in the partition contains at most one circled letter.  
\end{enumerate}
\begin{remark}
An $m\times n$ circled letter array may represent a preferred face of an $m\times n\times p$ proper array.  Hence, the $m\times n$ word representations form a proper subset of the $m\times n$ circled
letter arrays.
\end{remark}
\noindent {\bf Definition:}  Let $P_{m,n}$ be the number of $m\times n$ circled letter arrays.\\\\
\noindent {\bf Definition:}  Let $p_j(mn,c)$ be the $m\times n$ circled letter arrays that have $c$ letters and $j$ circles.\\\\
\noindent  Note that
\begin{equation}
P_{m,n} = \sum_{c=1}^{mn}\sum_{j=1}^{c}\left(\begin{array}{c}mn\\j\end{array}\right)p_j(mn,c)
\end{equation}
Furthermore, it is easy to show that 
\begin{equation}
p_j(mn,c) = \sum_{p=0}^{mn-c}\left(\begin{array}{c}mn-j\\p\end{array}\right)\left\{\begin{array}{c} n-j-p\\c-j\end{array}\right\}j^p
\end{equation}
where
\begin{itemize}
\item[1.] $p$ counts the number of times the circled letters reappear as unarrowed circled letters, i.e. the number of uncircled spaces contained in the union of partition elements that have a
cirlced letter.
\item[2.] $\left\{\begin{array}{c} n-j-p\\c-j\end{array}\right\}$ counts the number of ways to fill those spaces of the $m\times n$ grid that complement the union of partition elements containing
a circled letter.
 \end{itemize}
\noindent Remarks 2.2 and 2.3 provide useful algebraic representations for $\left\{\begin{array}{c} mn\\j\end{array}\right\}$ and $p_j(mn,c)$  [4],[7],[8].  
\begin{remark}
\[\sum_{mn=0}^{\infty}\left\{\begin{array}{c} mn\\t\end{array}\right\}\frac{y^{mn}}{(mn)!} = \frac{(e^y-1)^t}{t!}\]
\end{remark}
\begin{remark}
\[\sum_{mn-j}^{\infty}\frac{p_j(mn,c)z^{mn-j}}{(mn-j)!} = \frac{e^{zj}(e^z-1)^{c-j}}{(c-j)!}\]
\end{remark}
\begin{remark}
\noindent By using Remark 2.3 in Equation (2.1) and summing over $j$ and $c$, we can show the coefficient of $z^{mn}$ in the expansion of $\exp(e^z-1+ze^z)$ is $P_{m,n}$.
\end{remark}
\noindent With all the preliminary information in place, we are ready to develop a formula that provides an upper bound on the number of word representations
associated with $m\times n\times p$ proper arrays.  Define $W_{m,n}$ to be the number of
$m\times n$ circled letter arrays modulo $D_2$ symmetry.  Then, $W_{m,n}$ is our desired upper bound.  We calculate $W_{m,n}$ as follows.
\begin{itemize}
\item[1.] Let $S_{m,n}$ count the $m\times n$ circled letter arrays that are fixed via horizontal reflection, vertical reflection, and $180^{\circ}$ degree rotation.  
\item[2.] Let $H_{m,n}$ count the $m\times n$ circled letter arrays fixed via horizontal reflection.  Then, $H_{m,n} - S_{m,n}$ counts the $m\times n$ circled letter arrays that are fixed only by
horizontal reflection.  
\item[3.] Let $V_{m,n}$ count the $m\times n$ circled letter arrays fixed via vertical reflection.  Then, $V_{m,n} - S_{m,n}$ counts the $m\times n$ circled letter arrays that are fixed only by vertical
reflection. 
\item[4.] Let $R_{m,n}$ count the $m\times n$ circled letter arrays fixed via $180^{\circ}$ rotation.  Then, $R_{m,n} - S_{m,n}$ counts the $m\times n$ circled letter arrays that are fixed only by
rotation.  
\item[5.] Let $C_{m,n} = P_{m,n}-(H_{m,n}-S_{m,n})-(V_{m,n}-S_{m,n})-(R_{m,n}-S_{m,n})-S_{m,n}$.  Then, $C_{m,n}$ counts the $m\times n$ circled letter arrays that are not fixed by any symmetry
transformation.
\end{itemize}
\begin{theorem}
Let $W_{m,n}, C_{m,n}, P_{m,n}, V_{m,n}, H_{m,n}, R_{m,n},$ and $S_{m,n}$ be as previously defined.  Then 
\[ W_{m,n} =  \frac{P_{m,n}+H_{m,n}+V_{m,n}+R_{m,n}}{4}\]
\end{theorem}
{\bf Proof of Theorem 2.1:} To calculate the number of $m\times n\times p$ circled letter arrays modulo $D_2$ symmetry, we first determine whether a given $m\times n$ circled letter array,
called $A$, is fixed via any of the four symmetry transformations.  If $A$ is not fixed by any symmetry, it has four equivalent images.  However, if
$A$ is fixed under a symmetry transformation, it has at most two symmetry images.  It follows that  
\begin{align*}
W_{m,n} &= \frac{C_{m,n}}{4} + \frac{(H_{m,n}-S_{m,n})+(R_{m,n}-S_{m,n})+(V_{m,n}-S_{m,n})}{2} + S_{m,n}\\
&= \frac{P_{m,n}+H_{m,n}+V_{m,n}+R_{m,n}}{4} \qquad \Box
\end{align*}
\begin{remark}
Theorem 2.1 can be considered to be an immediate consequence of Burnside's Lemma.
\end{remark}
\subsection{A Numerical Example}
\noindent In order to understand how the information provided by Theorem 2.1 helps provide an upper bound on the basis size of the transition matrix, look at the following example.  Let $m = 3$ and $n =
1$.  The goal is to compute the size of the transition matrix associated with the $3\times 1\times p$ proper arrays.  It can be shown that the transition matrix is a $16\times 16$ matrix [1],[2].  Thus,
there are 16 possible $3\times 1$ {\it word} representations modulo symmetry.  Theorem 2.1 simply counts the number of symmetrically distinct $m\times n$ cirlced letter
arrays.  Ignoring symmetry, we calculate $P_{3,1} = 30$.  From additional calculations, we find that $H_{3,1} = 8,  V_{3,1} = 30,$ and  $R_{3,1} = 8$.  Theorem 2.1
implies that the number of $3\times 1$ circled letter arrays modulo symmetry is $\frac{30 + 8 + 30 + 8}{4} = 19$.  Thus, the transition matrix associated with the $3\times 1\times p$ proper arrays can
be no larger than a $19\times 19$ matrix.  The goal of our research is to obtain a formula that calculates the actual basis size of the transition matrix.  As this example demonstrates, the theorems in
this paper provide not actual basis size, but an upper bound.  
\section{Generating Function for Horizontal/Vertical Symmetry}
\noindent In order to use Theorem 2.1, we need to find generating functions for $H_{m,n}$, $V_{m,n}$, and $R_{m,n}$, These generating functions are obtained by dividing the
$m\times n$ array into small sections determined by the symmetry transformation.  For instance, in the case of horizontal reflection, the array is divided into two halves.  If
two or more transformations are applied, the array is divided into four quarters.  In either case, we can arbitrarily fill one of the halves/quarters with any arrangement of circled letters that obey
Remark 1.2, and then use symmetry to fill the remaining half/quarters of the $m\times n$ array.  The key in this technique is to carefully divide the array around the row and/or column that may
be fixed under a symmetry transformation.  Hence, the generating functions depend on the parity of
$m$ and $n$.\\\\
\noindent We begin finding the generating functions associated with $H_{m,n}$.  Due to symmetry, they are also the generating functions associated with $V_{m,n}$.  
We first assume the number of rows is even. No row or column is fixed by horizontal reflection.  We arbitrarily fill the
first $m$ rows with an arrangement of cirlced letters that obeys Remark 1.2.  Then, using horizontal reflection, we transfer the letters into the bottom $m$ rows.  If the partition element represented by
the letter does not contain a circle, there three possible ways to reflect that letter into the bottom $m$ rows.  Otherwise, since the partition element contains a circled letter, there are only two
possible ways to do this reflection.
\newpage
\noindent  Hence,
\[H_{2m,n} =
\sum_{j=1}^{mn}\sum_{q=0}^j\sum_{t=0}^{min(j-q,q)}\sum_{w=0}^{[\frac{j-q-t}{2}]}\frac{(mn)!p_q(mn,j)(j-q)!2^{j-q-t-3w}}{(mn-q)!t!(q-t)!w!(j-q-t-2w)!}\]
where,
\begin{enumerate}
\item[1.]  $j$ counts the letters in the first $m$ rows.  We will refer to the first $m$ rows as the {\bf top half} of the array and the remaining $m$ rows as the {\bf bottom half} of
the array.
\item[2.] $q$ counts the circles present among the $j$ letters.  
\item[3.] Let $A$ and $B$ be two letters in the top half of the array.  We say {\bf {\boldmath $A$} interchanges with {\boldmath $B$}} if $A$ maps to $B$ under a symmetry transformation.  For this
case, the symmetry transformation is horizontal reflection.  Then, $w$ counts the interchanges between uncircled letters.
\item[4.] $t$ counts the interchanges between a circled letter and an uncircled letter.
\end{enumerate}
\noindent Figure 3.1 illustrates the meaning of $j, q, w,$ and $t$.
\begin{center}
\vspace{-1.5cm}
\includegraphics[width=5.5cm, height=4.75cm]{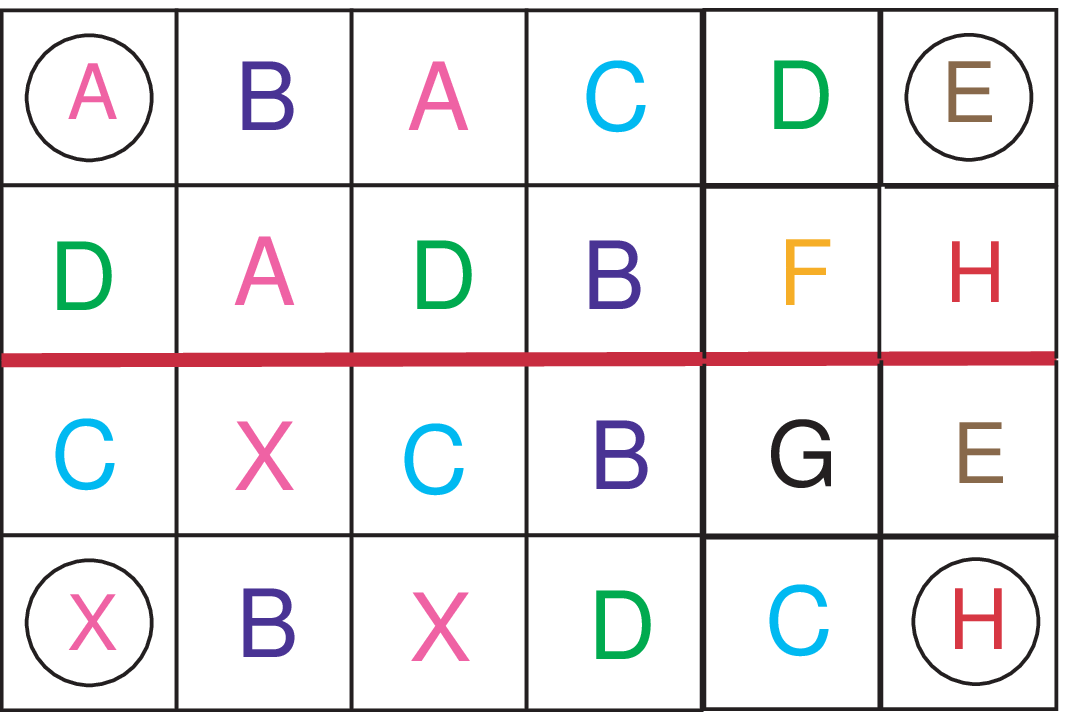}
\end{center}
\begin{center}
\small{Figure 3.1: This diagram is an example of a $4\times 6$ word representation that is fixed under horizontal reflection.  The top two rows are filled by seven letters, two of which are
circled.  We first discuss the three ways uncircled letters are transformed by horizontal reflection.  First, the $B$ reflects to itself.  Secondly, the $F$ reflects to $G$, where $G$ does not occur in
the top two rows.  Thirdly, we let $C$ reflect to $D$, where both $C$ and $D$ occur in the top two rows.   Next, we discuss the two ways circled letters are transformed by horizontal reflection. 
First, the $A$ can reflect to $X$.  Otherwise, the $H$ interchanges with $E$.}
\end{center}
\begin{theorem}
The coefficient of $z^{mn}$ in the expansion of $\exp(2(e^z-1)+\frac{(e^z-1)^2}{2}+ze^{2z})$ is $H_{2m,n}$.
\end{theorem}
\noindent {\bf Proof of Theorem 3.1:}  Let $k = mn$.  Then,
\[H(z) = \sum_{\scriptstyle k,j,q,t,w = 0}^{\infty}\frac{p_q(k,j)(j-q)!2^{j-q-t-3w}z^k}{(k-q)!t!(q-t)!w!(j-q-t-2w)!}\]
To obtain the desired result, use Remark 2.2 and sum over $j, q, t,$ and $w$.  See the proof of Theorem 3.2 for a more elaborate example.\qquad$\Box$\\\\
\noindent  The next sum, $H_{2m+1,n}$, counts the number of $(2m+1)\times n$ circled letter arrays obeying Remark 1.2 that are fixed under horizontal reflection.  Finding $H_{2m+1,n}$ is more complicated
since we must contend with a middle row that is fixed by the reflection map.  We obtain
 \[H_{2m+1,n}
=\sum_{j=1}^{mn}\sum_{i=0}^j\sum_{q=0}^{min(j-i,n)}\sum_{k=q}^{n}\sum_{E=0}^{q}\sum_{v=min(n-k,1)}^{n-k}\sum_{l=0}^{v}\sum_{t=0}^{min(j-i-q,i)}
\sum_{w=0}^{[\frac{j-i-q-t}{2}]}\]
\[\frac{(mn)!n!p_i(mn,j)p_E(k,q)p_l(n-k,v)2^{j-i-q-t-3w}(j-i)!}{(mn-i)!E!(k-E)!(n-k-l)!(i-t)!t!w!(j-i-q-t-2w)!l!}\]
\noindent where,
\begin{enumerate}
\setlength{\itemsep}{-\itemsep}
\item[1.]$j$ counts the letters in the top half of the array.
\item[2.] $i$ counts the circles that appear among the $j$ letters.  Hence, the top $m$ rows contain $j - i$ distinct uncircled letters.
\item[3.] $q$ counts how many of these $j-i$ uncircled letters occur in the middle row.
\item[4.] $k$ counts the positions in the middle row occupied by these $q$ letters.
\item[5.] $E$ counts the circles that occur among the $q$ letters.
\item[6.] $v$ counts the letters occuring in the remaining $n-k$ spaces of the middle row.   
\item[7.] $l$ counts the circles that occur among the $v$ letters.
\item[8.] $t$ counts the interchanges between the $i$ circled letters and the $j-i-q$ uncircled letters.
\item[9.] $w$ counts the interchanges between the $j-i-q-t$ uncircled letters.  
\end{enumerate}
\begin{center}
\includegraphics[width=6.0cm,height=3.5cm]{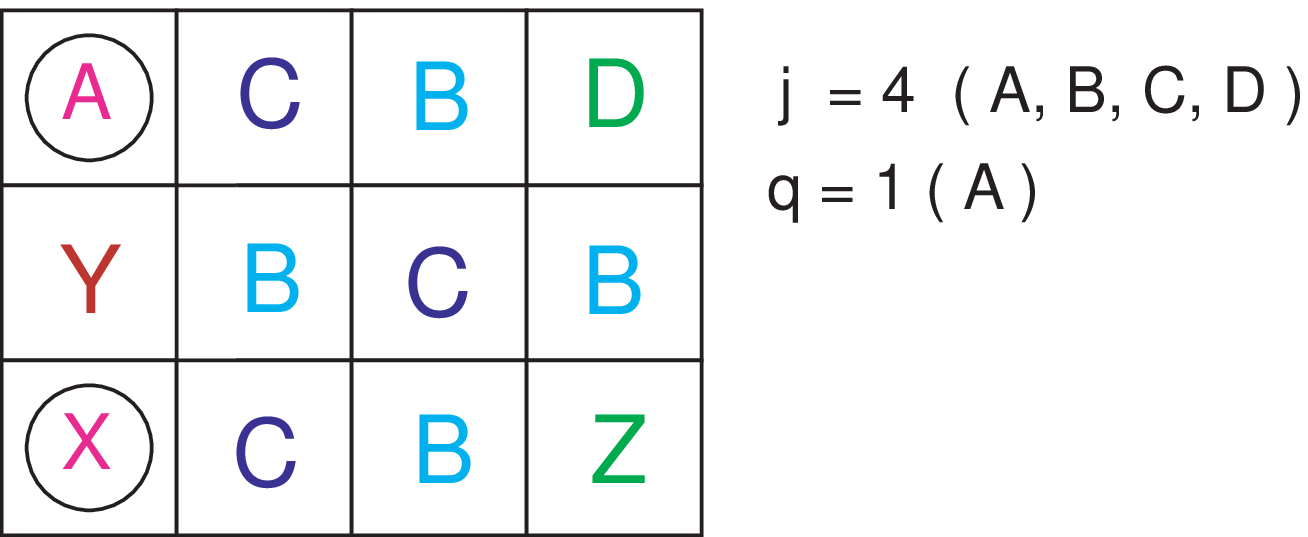}
\end{center}
\begin{center}
\small{Figure 3.2:  This diagram illustrates how uncircled letters that appear in the first $m$ rows may occur in the middle row.  For this particular example,
we see that $B$ and $C$ occur in both the top half of the array \underline{and} the middle row.  Therefore, $B$ and $C$ are fixed under horizontal reflection.} 
\end{center}
\newpage
\begin{theorem}
The coefficient of $z^{mn}y^n$ in the expansion of $\exp((y+1)e^{y+z}+(z+\frac{1}{2})e^{2z}-\frac{3}{2})$ is $H_{2m+1,n}$.
\end{theorem} 
\noindent{\bf Proof of Theorem 3.2:} Define
\[H(y,z) = \sum_{\scriptstyle mn,j,i,q,k,\atop \scriptstyle E,v,l,t,w,n = 0}^{\infty}\frac{p_i(mn,j)p_E(k,q)p_l(n-k,v)2^{j-i-q-t-3w}(j-i)!z^{mn}y^n}
{(mn-i)!E!(k-E)!(n-k-l)!(i-t)!t!w!(j-i-q-t-2w)!l!}\]
\noindent Using Remark 2.2, we find that
\begin{align*}
H(y,z) &= \sum_{\scriptstyle j,i,q, E, \atop \scriptstyle v,l,t,w = 0}^{\infty}
\frac{e^{zi}(e^z-1)^{j-i}e^{yl+yE}(e^y-1)^{v-l+q-E}z^iy^{E+l}2^{j-i-q-t-3w}}{(v-l)!(q-E)!E!(i-t)!t!w!(j-i-q-t-2w)!l!}\\
&= \sum_{\scriptstyle i,j,q, \atop \scriptstyle E,t,w = 0}^{\infty}
\frac{e^{ye^y+e^y-1}e^{zi}(e^z-1)^{j-i}e^{yE}(e^y-1)^{q-E}z^iy^E2^{j-i-q-t-3w}}{(q-E)!E!(i-t)!t!w!(j-i-q-t-2w)!}\\
&= \sum_{\scriptstyle i,q,E,t = 0}^{\infty}
\frac{e^{ye^y + e^y -1 +2(e^z-1) +\frac{1}{2}(e^z-1)^2}e^{zi}(e^z-1)^{q+t}e^{yE}(e^y-1)^{q-E}z^iy^E}{(q-E)!E!(i-t)!t!}\\
&= \sum_{\scriptstyle q,E = 0}^{\infty}
\frac{e^{ye^y + e^y - 1 + 2(e^z-1) + \frac{1}{2}(e^z-1)^2 +ze^{2z}}(e^z-1)^qe^{yE}(e^y-1)^{q-E}y^E}{(q-E)!E!}\\
&= e^{(y+1)e^{y+z}+(z+\frac{1}{2})e^{2z}-\frac{3}{2}} \qquad \Box
\end{align*}
\section{Generating Functions for Rotational Symmetry}
\noindent Our next step is to find a generating function for the rotational symmetry.  Observe that\\ $R_{2m,2n} = H_{2m,2n}$, $R_{2m+1,2n} = H_{2n,2m+1},$ and $R_{2m,2n+1} = H_{2m,2n+1}$.  It
remains to compute $R_{2m+1, 2n+1}$.  In this situation, the $180 ^{\circ}$ rotation fixes the central square.  Our geometric subdivision must avoid this central square.  Hence, the array is divided
into takes an upper $L$ and a lower $L$.
\begin{center}
\includegraphics[width=3.5cm, height=3.5cm]{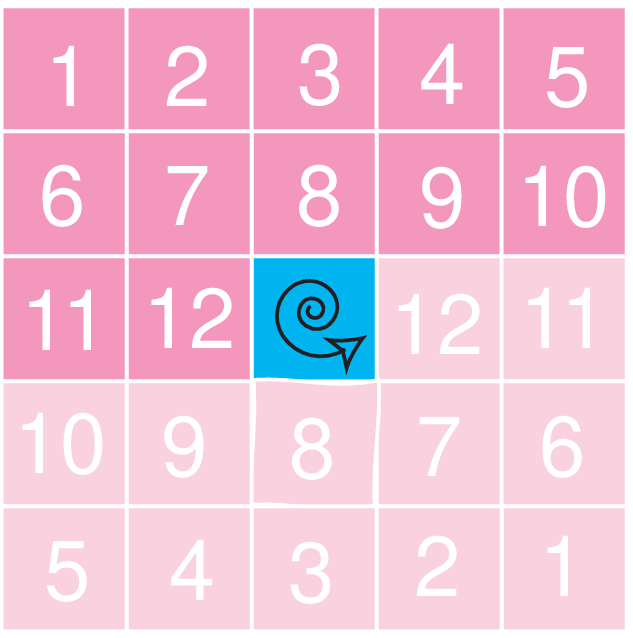}
\end{center}
\begin{center}
\small{Figure 4.1: Under $180^{\circ}$ rotation, the central square is fixed, and the array is subdivided into two $L$ shapes.}
\end{center}
\noindent We find that
\[R_{m,n}
=\sum_{j=1}^{N}\sum_{i=0}^{j}\sum_{q=0}^{j-i}\sum_{t=0}^{min(i,j-i-q)}\sum_{s=0}^{[\frac{j-i-q-t}{2}]}\frac{N!p_i(N,j)(j-i)!(1+q)}{(N-i)!t!q!(i-t)!s!2^{s-1}(j-i-q-t-2s)!}\]
\noindent where
\begin{enumerate}
\setlength{\itemsep}{-\itemsep}
\item[1.] $N =n\left[\frac{m}{2}\right] + \left[\frac{n}{2}\right]$
\item[2.] $j$ counts the letters that occur in the top $L$.
\item[3.] $i$ counts the circles that occur among these $j$ letters.  The top $L$ contains $j-i$ uncircled letters.  
\item[4.] $q$ counts the uncircled letters that rotate to themselves in the bottom $L$.  The central square can contain a letter that does not appear in the rest of the array, or it can be
contain one of these uncirlced letters enumerated by $q$. 
\item[5.] $t$ counts the interchanges between the $j-i-q$ uncircled letters and the $i$ circled letters.
\item[6.] $s$ counts the interchanges that occur among the $j-i-q-t$ .
\end{enumerate}
\begin{center}
t\includegraphics[width=10.0cm, height=3.0cm]{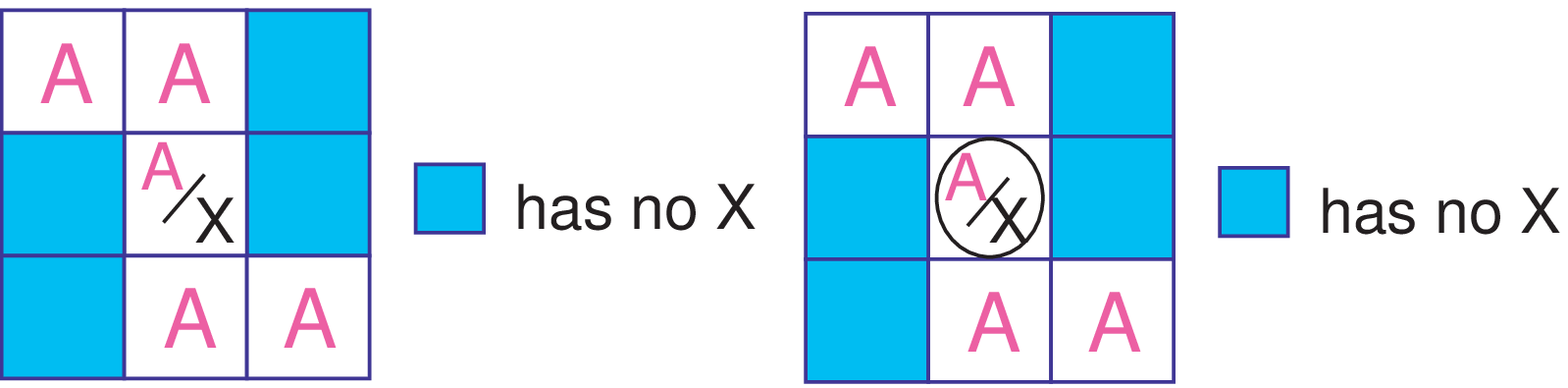}
\end{center}
\begin{center}
\small{Figure 4.2: The following diagram anaylzes the four ways the central square is completed.  In particular, it could be filled by an entirely new letter, $X$, or it could be filled by $A$, an
uncircled letter from the top $L$ that rotates to itself in the bottom $L$.  In either case, the middle position may be circled.}
\end{center}
\begin{theorem}
Let $m$ and $n$ be odd integers.  Let  $N = n[\frac{m}{2}] + [\frac{n}{2}]$.  The coefficient of $z^N$ in the expansion of $\exp(\frac{1}{2}(e^z-1)^2+2(e^z-1)+z+ze^{2z})$ is $R_{m,n}$
\end{theorem}
\noindent {\bf Proof of Theorem 4.1:} Let
\[R(z) =\sum_{\scriptstyle j,i,q, \atop \scriptstyle N,t,s = 0}^{\infty}\frac{p_i(N,j)(j-i)!(1+q)z^N}{(N-i)!t!q!(i-t)!s!2^{s-1}(j-i-q-t-2s)!}\]
\noindent By using the techniques of Theorem 3.1, we obtain the desired result. \qquad $\Box$.
\section{Generating Function for Self Symmetrical Word Representations}
\noindent We define an $m\times n$ word representation to be {\bf fully symmetric} if and only if it is fixed via horizontal \underline{and} vertical reflection.
 Let $S_{m,n}$ count the number of fully symmetrical word
representations.  Although $S_{m,n}$ only appears in the proof of Theorem 2.1, we can apply the techniques of Sections 3 and 4 to determine the generating functions for $S_{m,n}$.  Since the
fully symmetrical word representations are invariant under two symmetry transformations, the underlying $m\times n$ array is divided into four quadrants.  The upper left quadrant is arbitrarily
filled with an arrangement of circled letters which obey Remark 1.2.  Then, the remaining three quadrants are filled by applying the two symmetry transformations to the upper left quadrant.  There are two
types of letters in the upper left quadrant.  The first type of letter is called a  singleton letter.  A {\bf singleton  letter} is a letter whose image in the remaining three quadrants is never another
letter that appears in the upper left quadrant.  There are five possible ways an uncircled letter may be a singleton, and one way in which a circled letter may be a singleton.  These six
cases are illustrated in Figure 5.1.
\begin{center}
\includegraphics[width=6.0cm, height=4.0cm]{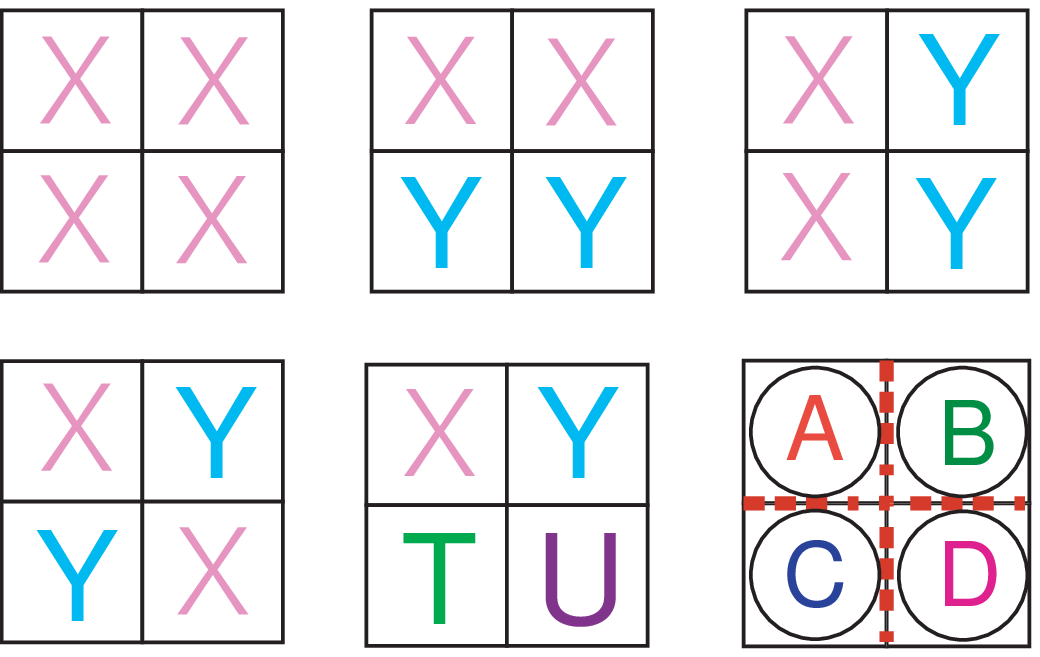}
\end{center}
\begin{center}
\small{Figure 5.1: The five ways an uncircled singleton letter can be transformed under symmetry and the unique way a circled letter is a singleton.}
\end{center}
\noindent  The second type of letter present in the upper left quadrant is considered to be part of a {\bf double pair}.  A letter is part of {\bf double pair}
when its image in one of the other three quadrants is another letter orginally present in the upper left quadrant.  Figure 5.2 illustrates the six ways a double
pair composed of two uncircled letters tranforms in a self symmetrical manner while Figure 5.3 illustrates the three possiblities for a double pair composed of one circled letter and one uncircled letter.
\begin{center}
\includegraphics[width=6.0cm, height=10.0cm]{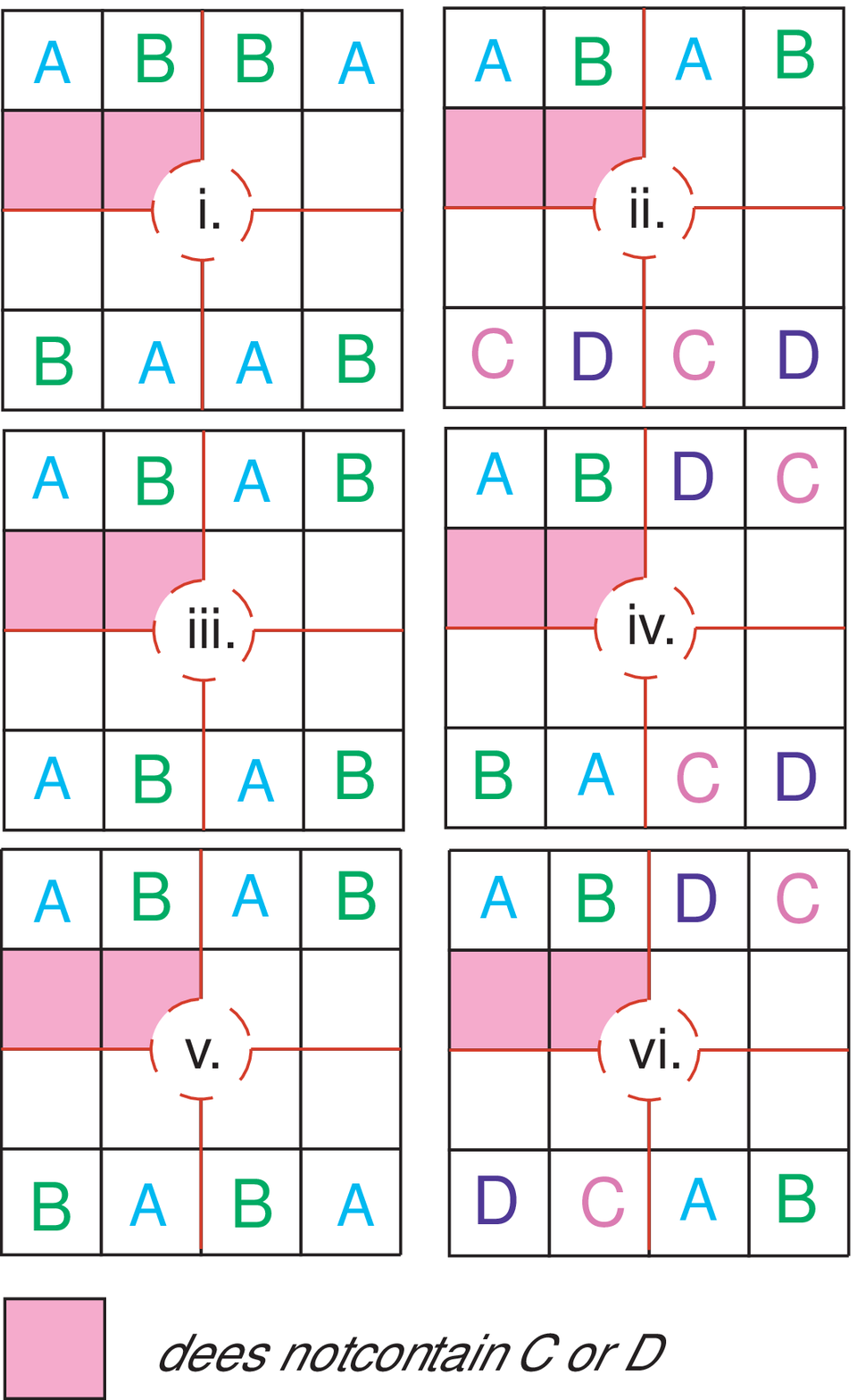}
\end{center}
\begin{center}
\small{Figure 5.2: The six ways two uncircled letters, each appearing in the upper left quadrant, may be interchanged during
symmetry operations.} 
\end{center}
\begin{center}
\includegraphics[width=8.0cm, height=3.5cm]{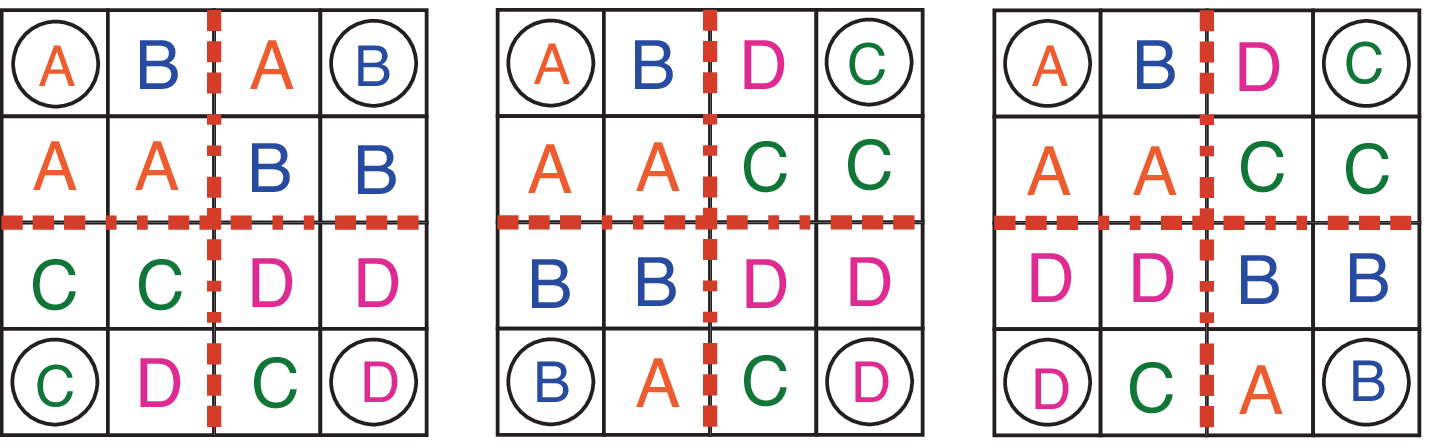}
\end{center}
\begin{center}
\small{Figure 5.3: Three ways the circled $A$ interchanges with the uncircled $B$.  Letters $C$ and $D$ do not appear in the
upper left quadrant but are images of $A$ and $B$ under horizontal and vertical reflection.}
\end{center}
\noindent We begin by assuming that the row and
column dimensions of the underlying array are both even.  Hence, there are no central rows or columns to worry about.  The geometric construction depends on locating the singleton letters and the double
pairs present in the upper left quadrant.  In particular,  
\newpage
\[S(2m,2n) =
\sum_{j=1}^{mn}\sum_{q=1}^j\sum_{l=0}^{min(q,j-q)}\sum_{s=1}^{[\frac{j-q-l}{2}]}\frac{(mn)!p_q(mn,j)3^{l+s}(j-q)!5^{j-q-l-2s}}{(mn-q)!l!(q-l)!s!(j-q-l-2s)!}\]
\noindent where
\begin{enumerate}
\setlength{\itemsep}{-\itemsep}
\item[1.] $j$ counts the letters in the upper left quadrant
\item[2.] $q$ counts the circles that occur among the $j$ letters.  There are $j-q$ uncircled letters that occur in the upper left quadrant.
\item[3.] $l$ counts the interchanges between the $q$ circled letters and the $j-q$ uncircled letters.
\item[4.] $s$ counts the interchanges between the $j-q-l$ circled letters.
\end{enumerate}
\noindent The proof of Theorem 5.1 is similiar to that of Theorem 3.1 and will be omitted.
\begin{theorem}
The coefficient of $z^{mn}$ in the expansion of $\exp(3(1+z)e^{2z}-(2z+1)e^z-2)$ is $S_{2m,2n}$.
\end{theorem}
\noindent When one of the row or column dimensions is odd, we must deal with a central column and/or row that is fixed by one of the symmetry transformations.  For word representations of dimensions
$2m\times (2n+1)$, we decompose the underlying array into four quadrants and a middle column.  Define the {\bf top half of the middle column} to be the first $m$ squares.  In order to compute
$S_{2m,2n+1}$, we begin by filling the top half of the middle column with an arbitrary arrangement of circled letters that obey Remark 1.2.  Notice that the letters in the top half of the middle column,
when reflected to the bottom half of this column, may either go to themselves, go to an entirely new image, or be part of a double pair.  In all three cases, letters present in the middle may fill the
four quadrants in a unique manner, modulo order.
\begin{center} 
\includegraphics[width = 12.0cm, height = 12.0cm]{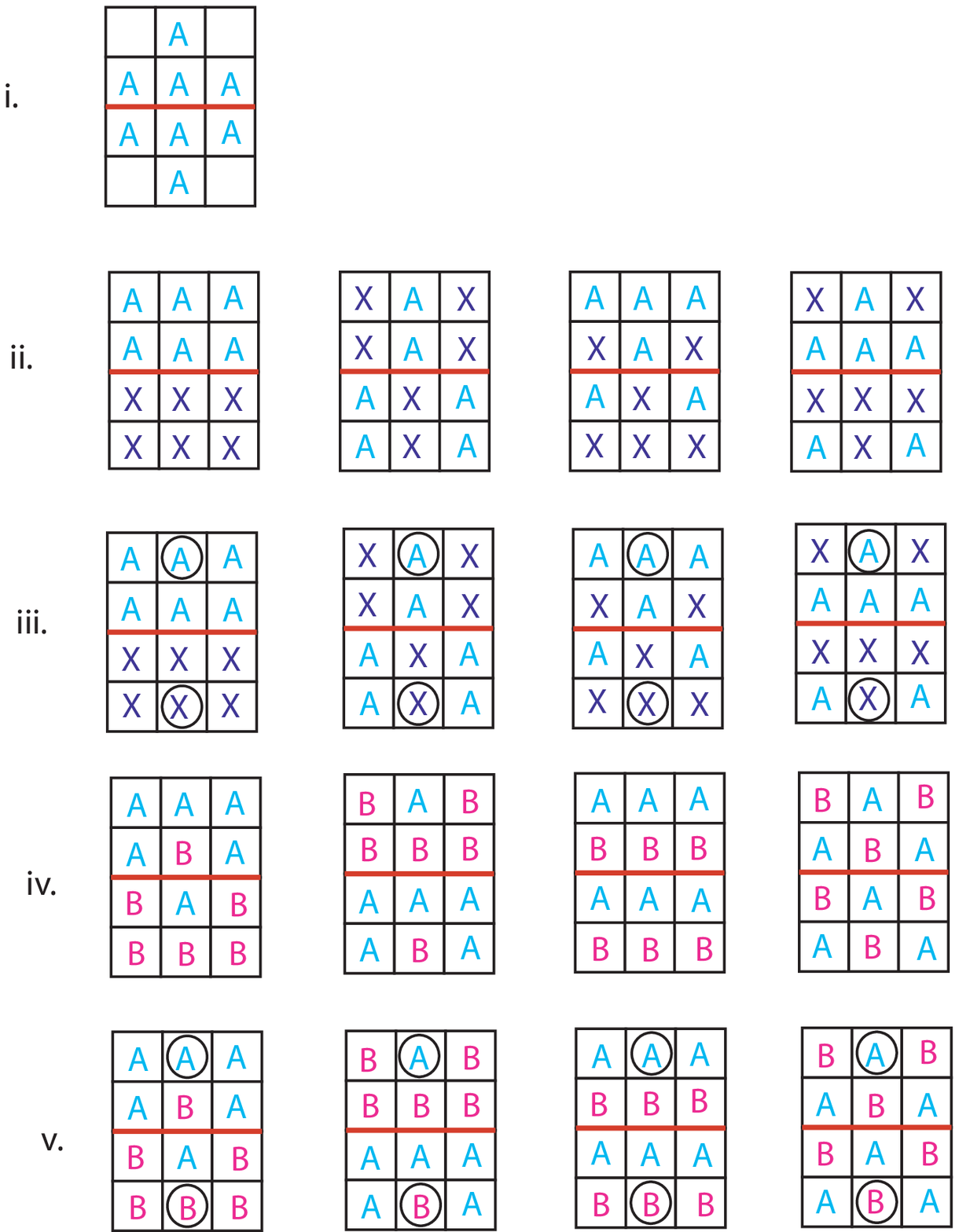}
\end{center}
\begin{center} \small{Figure 5.4: These diagrams illustrates how letters that appear in the middle column may appear in the quadrants.  We begin by analyzing how the letters in the top half of the
middle column reflect to the bottom half of that column.  In Diagram i, the letter in the top half of the middle column goes to itself under horizontal reflection.  The
letter may occur in the quadrants as a singleton.  In Diagrams ii and iii, the letter in the top half of the middle column, $A$, reflects to a new letter, $X$.  In
either case, both $A$ and its image, $X$, may appear in the quadrants as either a singleton or a double pair.  In Diagrams iv and v, the two letters present in the top half of the middle column
interchange positions under horizontal reflection.  Just like in Diagrams ii and iii, these two letters may occur in the quadrants as either a singleton or a double pair.}
\end{center}
\noindent  Once the letters that occur in the middle column are placed in the quadrants, we fill the remaining squares in the upper left quadrant with an arbitrary arrangement of circled letters that
obey Remark 1.2.  Then, we complete the remaining three quadrants in the manner illustrated by Figures 5.1 through 5.3.  Following the steps outlined above, we arrive at a summation of 19 variables. 
The derivation of the exponential generating function given in Theorem 5.2 utilizes the techniques of Theorom 3.2.  Details are available, upon request, from the author. 
\begin{theorem}
The coefficient of $z^mx^{mn}$ in the expansion of\\ $\exp((z+\frac{1}{2})e^{2x+2z}+(3x+\frac{5}{2})e^{2x}-2(x+1)e^x+e^{z+x}-2)$ is $S_{2m,2n+1}$.
\end{theorem}
\noindent The final sum of this section counts the number of $(2m+1)\times (2n+1)$ fully symmetrical word representations.  In the case the geometric decomposition involves four quadrants, a central
square, and a middle cross, where a {\bf middle cross} is the union of the middle column and middle row minus the central square.   Thus, a middle cross consists of the top and bottom halves of the middle
column along with the left and right halves of the middle row.  In this situation, {\bf the top half of the middle column} is the $m$ squares above the central square, and {\bf the left half of the middle
row} is the $n$ squares to the left of the central square.\\\\ 
\noindent Fortunately, the strategy for enumerating the $(2m+1)\times (2n+1)$ fully symmetrical word representations is an adaptation of the techniques we used for determining $S_{2m,
2n+1}$.  Hence, we will simply outline the procedure for calculating $S_{2m+1, 2n+1}$ and then provide the resulting generating function.  The major difference between the two cases is the fixed middle
row with its central square.  We once again start by filling in the top half of the middle column.  Using horizontal reflection, we are able to determine the bottom half of this middle column. 
Next, we complete the left half of the middle row.  Note that if a letter from the top half of the middle column has the same image in the bottom half of the middle column, it
\underline{may} occur in the middle row.  Using vertical reflection, we are then able to fill in the right half of the middle column in a manner similiar to the one used in the completion of the middle
row.   To visualize how this is done,  simply rotate Figure 5.4 by 90 degrees.  Once the middle column and row are completed, we
turn to the quadrants.  Letters that appear in the middle cross may occur, modulo order, in a unique fashion in the quadrants.  Once the letters from the middle cross are transferred to
the quadrants, we fill the remaining spaces in the quadrants via the technique discussed in the construction of $S_{2m,2n}$.  Finally, we deal with the fixed central square.  This position may be filled
with a singleton or an entirely new letter.  In either case, this letter may be circled.
\begin{center}
\includegraphics[width=5.0cm, height=5.0cm]{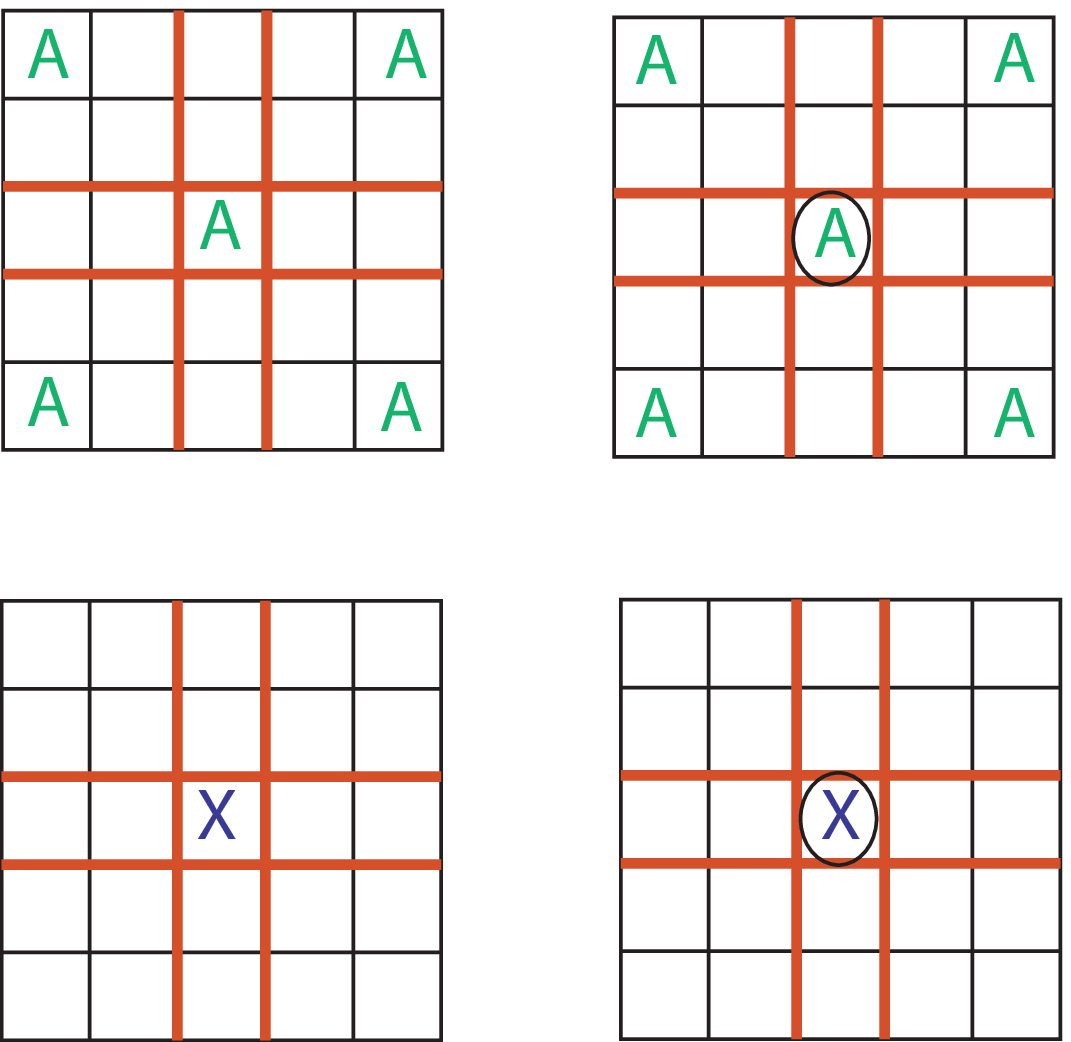}
\end{center}
\begin{center}
\small{Figure 5.5:  The four ways the middle position may be filled in the case of a fully symmetric word representation.  Note that $X$ does not appear any place else but in
the middle position while $A$ is fixed under both horizontal and vertical reflection.}
\end{center}
\noindent Following the steps outlined above, we arrive at a summation of 35 variables.  The derivation of the exponential generating function utilizes the previously mentioned techniques.  Details are
available, upon request, from the author.  The results of these calculations are recorded in Theorme 5.3.  
\begin{theorem}
The coefficient of $x^my^nz^{mn}$ in the expansion of\\ 
$\exp(y + z + x - 2 + ye^{2y + 2z} + xe^{2x + 2z} + 3ze^{2z} - 2ze^z + e^{x + y + z} - 2e^z + 2e^{2z} +\frac{1}{2}e^{2x + 2z} +\frac{1}{2}e^{2y + 2z})$ is $S_{2m+1,2n+1}$.
\end{theorem}
\section{Open Questions}
\noindent By using a particular decomposition of the Bell Numbers [4],[7],[8] and applying various symmetry transformations to $m\times n$ rectangular arrays,
we discovered a formula that calculates an upper bound for the number of $m\times n$ word representations.  In order to exactly predict the number of $m\times n$ word representations, we need to apply an
adjacency rule [1],[2] to the set of $m\times n$ arrays of circled letters enumerated by Theorem 2.1.  Thus, future research should concentrate on how this adjacency condition reduces the upper bound
provided by Theorem 2.1.\\\\  
\noindent Another promising avenue of research involves square arrays, where $m = n$.  Since a square array has $D_4$ symmetry, the formula to enumerate $m\times m$ letter representations modulo symmetry
will be more complicated than the corresponding formula provided by Theorem 2.1.  We are in the process of calculating the
generating functions necessary in the enumeration of $m\times m$ arrays of circled letters.  Once this is
accomplished, an upper bound for the basis size of the associated transition matrix will be completely determined for all integers $m$ and $n$.\\\\
\noindent Finally, a different type of research would involve exploring the connections between $m\times n\times p$ proper arrays and percolation theory.  At the present time, we have not explored
the connection in any depth but realize that the stochastic and probabilisitic techniques of percolation theory could, when applied to the representation of an $m\times n\times p$ proper array as a bond
percolation on $Z^3$ with an open cluster at the origin (see Section 0.2), give rise to a whole new category of results.
\section*{Acknowledgments}
\noindent The author would like to thank Dr. Bruce Sagan and Dr. Robert Sulanke for their suggestions regarding the Introduction.  The author also thanks Dr. Harris Kwong for his time and help in editing
the previous drafts of this paper.
\newpage
\section*{Appendix A: Numerical Data}
\noindent The following table provides, for small integer values of $m$ and $n$, numerical values for $P_{m,n}, H_{m,n}, V_{m,n},$\\
$R_{m,n},$, and $S_{m,n}$. All the values
came from the generating functions given by Theorems 3.1, 3.2, 4.1, 5.1, 5.2, and 5.3 and were verified by a short Maple program the author created.
\begin{center}
\begin{tabular}{||l|c|c|c|c|c|}
\hline
$m\times n$ & $P_{m,n}$ & $H_{m,n}$ & $V_{m,n}$ & $R_{m,n}$ & $S_{m,n}$\\ \hline
$2\times 3$ & 5653 & 107 & 197 & 107 & 23\\ \hline
$2\times 4$ & 306419 & 851 & 851 & 851 & 55\\ \hline
$2\times 5$ & 22277080 & 770 & 12976 & 770 & 234 \\ \hline
$3\times 2$ & 5653 & 197 & 107 & 107 & 23 \\ \hline
$3\times 4$ & 2062199125 & 463973 & 79525 & 79525 & 1525\\ \hline
$3\times 5$ & 2678973711602 & 35802956 & 8315630 & 3302472 & 26168 \\ \hline
\end{tabular}
\end{center}
\begin{center}
\small{Table 2:  Numerical Data for certain $m\times n$ word representations}
\end{center}
\newpage
\section*{References}
$[1]$\hspace{0.2cm} J. Quaintance ``$m\times n$ Proper Arrays: Geometric and Algebraic Methods of Classification'', Ph.D. Dissertation, University of Pittsburgh, August 2002\\\\
$[2]$\hspace{0.2cm} J. Quaintance ``$n\times m$ Proper Arrays: Geometric Constructions and the Associated Linear Cellular Automata'', {\it Maple Summer Workshop 2004 Procedings}\\\\
$[3]$\hspace{0.2cm} J. Quaintance ``Combinatoric Enumeration of the Geometric Classes Associated with $n\times p$ Proper Arrays'', preprint January 2004\\\\ 
$[4]$\hspace{0.2cm} N. J. A. Sloane, ``The On-Line Encyclopedia of Integer Sequences'',\\ http://www.research.att.com/~njas/sequences/ \\\\
$[5]$\hspace{0.2cm} S. Lo and M. Monagan, `` A Modular Algorithm for Computing the Characterisitic Polynomial of an Integer Matrix in Maple '', {\it Maple Summer Workshop 2005 Procedings}\\\\
$[6]$ \hspace{0.2cm} J. Quaintance `` Letter Representations of Rectangular $m\times n\times p$ Proper Arrays", preprint June 2004\\\\
$[7]$\hspace{0.2cm} K. Yoshinaga and M. Mori, `` Note on an Exponential Generating Function of Bell Numbers '', {\it Bull. Kyushu Inst. Tech.}, {\bf 24}, 23-27\\\\
$[8]$\hspace{0.2cm} D. Branson, ``Stirling Numbers and Bell Numbers: Their Role in Combinatorics and Probability '', {\it Math. Scientist}, {\bf 25}, 1-31\\\\
$[9]$\hspace{0.2cm} B. Cipra, `` An Introduction to the Ising Model '', {\it The American Mathematical Monthly}, {\bf 94}, No. 10, 937-959\\\\
$[10]$\hspace{0.2cm} G. Grimmett, {\it Pecolation} Second Edition, Springer-Verlag, 1991\\\\
$[11]$\hspace{0.2cm} D. Stauffer and A. Aharony, {\it Introduction to Percolation Theory} Second Edition, Taylor and Francis, 1992\\\\
$[12]$\hspace{0.2cm} R. Stanley, {\it Enumerative Combinatorics, Volume 1} New York: Wadsworth and Brook/Cole, 1986 \\\\
\newpage
\end{document}